\RequirePackage{ifpdf}
\ifpdf 
\documentclass[pdftex]{sigma}
\else
\documentclass{sigma}
\fi

\numberwithin{equation}{section}
\newcommand{\bs}{\boldsymbol}

\renewcommand{\a}{\alpha}
\renewcommand{\b}{\beta}
\renewcommand{\L}{\Lambda}
\renewcommand{\l}{\lambda}
\renewcommand{\d}{\delta}
\newcommand{\D}{\Delta}
\newcommand{\G}{\Gamma}
\newcommand{\g}{\gamma}

\newcommand{\bb}{{\mathbf b}}
\newcommand{\boe}{{\mathbf e}}
\newcommand{\bi}{{\mathbf i}}
\newcommand{\bj}{{\mathbf j}}
\newcommand{\bk}{{\mathbf k}}
\newcommand{\bm}{{\mathbf m}}
\newcommand{\br}{{\mathbf r}}
\newcommand{\bu}{{\mathbf u}}
\newcommand{\bx}{{\mathbf x}}
\newcommand{\boy}{{\mathbf y}}

\renewcommand{\Re}{\mathop{\mathrm{Re}}\nolimits}
\newcommand{\Det}{\mathop{\mathrm{Det}}\nolimits}

\begin{document}

\allowdisplaybreaks

\renewcommand{\PaperNumber}{119}

\FirstPageHeading

\ShortArticleName{A System of Multivariable Krawtchouk Polynomials and a Probabilistic Application}

\ArticleName{A System of Multivariable Krawtchouk Polynomials\\ and a Probabilistic Application}

\Author{F. Alberto GR\"UNBAUM~$^\dag$ and Mizan RAHMAN~$^\ddag$}

\AuthorNameForHeading{F.A.~Gr\"unbaum and M.~Rahman}

\Address{$^\dag$~Department of Mathematics, University of California, Berkeley, CA  94720, USA}
\EmailD{\href{mailto:grunbaum@math.berkeley.edu}{grunbaum@math.berkeley.edu}}
\URLaddressD{\url{http://www.math.berkeley.edu/~grunbaum/}}

\Address{$^\ddag$~Department of Mathematics and Statistics, Carleton University,\\
\hphantom{$^\ddag$}~Ottawa, ONT, Canada, K1S~5B6}
\EmailD{\href{mailto:mrahman@math.carleton.ca}{mrahman@math.carleton.ca}}

\ArticleDates{Received June 10, 2011, in f\/inal form December 19, 2011; Published online December 27, 2011}

\Abstract{The one variable Krawtchouk polynomials, a special case of the $_2F_1$ function
did appear in the spectral representation of the transition kernel for a~Markov chain studied a~long time ago by M.~Hoare and M.~Rahman. A multivariable extension of this Markov chain was considered in a later paper by these authors where a certain two variable extension of the $F_1$ Appel function shows up in the spectral analysis of the corresponding
transition kernel. Independently of any probabilistic consideration a certain
multivariable version of the Gelfand--Aomoto hypergeometric function was considered in papers by H.~Mizukawa and H.~Tanaka. These authors and others such as P.~Iliev and P.~Tertwilliger treat the two-dimensional version of the Hoare--Rahman work  from a Lie-theoretic point of view. P.~Iliev then treats the general $n$-dimensional case. All of these authors proved several properties of these functions. Here we show that these functions play a crucial role
in the spectral analysis of the transition kernel that comes from pushing the work of Hoare--Rahman to the multivariable case. The methods employed here to prove this as well as several properties of these functions are completely dif\/ferent to those used by the authors mentioned above.}


\Classification{33C45; 22E46; 33C45; 60J35; 60J05}

\vspace{-2mm}

\section{Introduction}
\label{sec1}

\looseness=-1
The genesis of this paper goes back to a joint work of Hoare and Rahman \cite{HR1} in 1983, where the idea of the so-called ``Cumulative Bernoulli Trials'' (CBT) was introduced.  The essential elements of this probabilistic model are as follows:  A player (say, of poker dice) rolls a subset,~$i$, of a f\/ixed number, $N$, of dice for success (say, ``aces'') with a certain probability $\a$.  The player is allowed to save his/her $k$ successes, and given a second chance, namely, to mix the $i-k$ unsuccessful dice with the previously $N-i$ unrolled ones.
 The player then rolls the combined dice, numbering $N-k$ for success with probability~$\b$.
Suppose the number of successes in this second try is $j-k$, which, when combined with the previously earned points, $k$, gives the total number of successes as $j$. If the number of successes after these two rolls is def\/ined to be the state of our system we get a Markov chain by iterating this scheme which took as from state~$i$ to state~$j.$

More explicitly the transition probability matrix of this Markov chain with state space ${0,1,2,\dots,N}$ is given by
\begin{gather}
\label{eq1.1}
K(j,i) = \sum_{k=0}^{\min(i,j)} b(k;i;\a)b(j-k;N-k;\b),
\end{gather}
where
\begin{gather}
\label{eq1.2}
b(k;n;p) = \binom{n}{k} p^k(1-p)^{n-k},
\end{gather}
is the binomial distribution.  The stationary distribution $\phi_0(i)$ corresponding to this process can be def\/ined by
\begin{gather*}
\sum_{j=0}^N K(i,j)\phi_0(j) = \phi_0(i).
\end{gather*}

A suf\/f\/icient condition for some $\phi_0(i)$ to satisfy this condition is
\begin{gather}
\label{eq1.4}
K(i,j)\phi_0(j) = K(j,i)\phi_0(i),
\end{gather}
and also that the summation part of $K(j,i)$ is symmetric in $i$ and $j$.  Use of \eqref{eq1.1}, \eqref{eq1.2} and~\eqref{eq1.4} gives
\begin{gather*}
\phi_0(i) = b(i;N;\eta),
\end{gather*}
where
\begin{gather*}
\frac {(1-\a)\eta}{\b} = \frac {1-\eta}{1-\b} = D_1^{-1},
\end{gather*}
with $D_1 = 1 + \frac {\a\b}{1-\a}$.

As determined in \cite{HR1} the eigenvalues $\l_k$ of the eigenvalue equation
\begin{gather*}
\sum_{j=0}^N K(i,j)\psi_k(j) = \l_k\psi_k(i),\qquad k = 0,1,2,\dots,
\end{gather*}
are $\l_k = \a^k(1-\b)^k$, with eigenfunctions
\begin{gather*}
\psi_k(i) = b(i;N;\eta)\, {}_2F_1\big({-}i,-k;-N;\eta^{-1}\big),
\end{gather*}
where the hypergeometric function ${}_2F_1$ in one variable $i$ is just the Krawtchouk polynomials, see~\cite{HTF}.  Clearly, by use of the orthogonality property of Krawtchouk polynomials we can write down the spectral representation of $K(j,i)$, namely,
\begin{gather*}
K(j,i) = b(j;N;\eta) \sum_{k=0}^N \binom{N}{k}^{-1} \left( \frac {(1-\a)(1-\b)}{\b}\right)^k (\a(1-\b))^k \nonumber\\
\phantom{K(j,i) =}{}\times {}_2F_1\big({-}i,-k;-N;\eta^{-1}\big)\, {}_2F_1\big({-}j,-k;-N;\eta^{-1}\big).
\end{gather*}

It is obvious that the above prototype of dice-tossing and ``saving'' successes can be extended, on one hand, to more practical situations in which Bernoulli trials may be accumulated (for example, in `infection-therapy' models, see \cite{HR1,HR2}), and on the other hand, to multiple variables where one def\/ines various kinds of success (say, aces, kings, queens, $\dots$).  Suppose we have a~process where the total number of dice is $N$, of which $n$ subsets $i_1,\dots,i_n$ are tossed separately for successes of $n$ dif\/ferent kinds with probabilities $\a_1,\a_2,\dots,\a_n$, respectively.  Let $k_1,k_2,\dots,k_n$ be the number of successes in each category (that is, $k_1$ aces, $k_2$ kings, $\dots$, etc.).  Mix the ``unsuccessful'' $i_1-k_1,i_2-k_2,\dots,i_n-k_n$ trials with the remaining $N-i_1-i_2-\dots-i_n$ dice.  The player ``saves'' the $\bk$ successes and is allowed to try again with the resulting $N-k_1-\dots-k_n$ dice with probabilities $\b_1,\b_2,\dots,\b_n$ of producing $j_1-k_1$ aces, $j_2-k_2$ kings, \dots, etc.

With the def\/inition for the $n$-fold multinomial distribution $b_n$ given below we get
 that the transition probability kernel from the state $(i_1,i_2,\dots,i_n)$ to $(j_1,j_2,\dots,j_n)$ is
\begin{gather}
 K(j_1,\dots,j_n;i_1,\dots,i_n)  \label{eq1.10}\\
{}
 = \sum_{k_1} \dots \sum_{k_n} b_n(j_1-k_1,\dots,j_n-k_n;N-k_1-k_2-\dots-k_n;\b_1,\dots,\b_n)  \prod_{r=1}^n b(k_r;i_r;\a_r).\nonumber
\end{gather}

The stationary distribution in this case, as before, is def\/ined by
\begin{gather}
K(\bj;\bi)\phi_0(\bi) = K(\bi;\bj)\phi_0(\bj), \qquad \bi = (i_1,i_2,\dots,i_n),\qquad \bj = (j_1,j_2,\dots,j_n),\label{eq1.11}
\end{gather}
which, together with \eqref{eq1.10} gives, as a suf\/f\/icient condition
\begin{gather*}
\phi_0(\bi)  = b_n(\bi;N;\bs{\eta})
 \equiv \binom{N}{i_1,\dots,i_n} \prod_{k=1}^n \eta_k^{i_k} \left( 1 - \sum_{j=1}^k \eta_j\right)^{N - \sum\limits_{j=1}^k i_j},
\end{gather*}
which is the $n$-fold multinomial distribution. Using \eqref{eq1.10} and \eqref{eq1.11} we f\/ind that, in perfect analogy to the one-variable case, the $\eta$'s are related to the probability parameters $\a$'s and $\b$'s in the following way:
\begin{gather}
\frac {1-\a_1}{\b_1} \eta_1  = \frac {1-\a_2}{\b_2} \eta_2 = \cdots = \frac {1-\a_n}{\b_n} \eta_n
 = \frac {1-\overset{n}{\underset{k=1}{\sum}}\eta_k}{1-\overset{n}{\underset{k=1}{\sum}}\b_k} = 1 - \sum_{k=1}^n \a_k\eta_k = D_n^{-1},\label{eq1.13}
\end{gather}
where
\begin{gather*}
D_n = 1 + \sum_{k=1}^n \frac {\a_k\b_k}{1-\a_k},\qquad 0 < \a_k,\, \b_k < 1.
\end{gather*}

One of the questions we will address in this paper is this:  what are the eigenvalues and eigenfunctions of $K(\bi;\bj)$?  In other words, we will seek solutions of the eigenvalue problem:
\begin{gather}
\label{eq1.15}
\sum_{j_1,\dots,j_n=0}^N K(\bi;\bj)\psi_{\bk}(\bj) = \l_{\bk}\psi_{\bk}(\bi).
\end{gather}

In the single-variable case the eigenfunctions are simply $\phi_0(\bi)$ times the ordinary Krawtchouk polynomials, as we mentioned earlier.  So it is reasonable to expect that in $n$ dimensions $(n \ge 2)$ the solutions will be an appropriate extension of these polynomials. The question is: which one?  Even in the $n = 2$ case the question is not quite as straightforward as it would seem.  In fact, the same authors, Hoare and Rahman, wrestled with this problem for a number of years until they were able to show that the Krawtchouk limit of the $9-j$ symbols of quantum angular momentum theory in physics provides the answer which, written in slightly more convenient notation is
\begin{gather}
\label{eq1.16}
b_2(x_1,x_2;N;\eta_1,\eta_2)F_1^{(2)}(-m_1,-m_2;-x_1,-x_2;-N;t,u,v,w),
\end{gather}
where $(t,u,v,w)$ satisfy certain relationships with $\eta_1$ and $\eta_2$, and the $F_1^{(2)}$ represents an iterate or a $2$-variable extension of the familiar $F_1$ Appell function:
\begin{gather*}
F_1(a;b,b';c;x,y) = \sum_i \sum_j \frac {(a)_{i+j}(b)_i(b')_j}{i!j!(c)_{i+j}} x^iy^j,
\end{gather*}
that is,
\begin{gather*}
F_1^{(2)}(a,a';b,b';c;t,u,v,w) = \sum_k\sum_j\sum_k\sum_l \frac {(a)_{i+j}(a')_{k+l}(b)_{i+k}(b')_{j+l}}{i!j!k!l!(c)_{i+j+k+l}} t^iu^jv^kw^l,
\end{gather*}
subject to conditions for convergence in case they are inf\/inite sums.  Unbeknownst to the authors, Hoare and Rahman, at the time they published the paper \cite{HR2}, a general $n$ variable extension of the $F_1^{(2)}$ in~\eqref{eq1.16}, namely,
\begin{gather}
\label{eq1.19}
F_1^{(n)}(-\bm,-\bx;-N;\bu) := \sum_{\underset{i,j}{\sum} k_{ij} \le N} \frac {\overset{n}{\underset{i=1}{\prod}} (-m_i)_{\overset{n}{\underset{j=1}{\sum}} k_{ij}} \overset{n}{\underset{i=1}{\prod}} (-x_i)_{\overset{n}{\underset{j=1}{\sum} } k_{ji}}}{\underset{i,j}{\prod} k_{ij}!(-N)_{\underset{i,j}{\sum} k_{ij}}}  \prod_{i,j} u_{ij}^{k_{ij}},
\end{gather}
a special case of Gelfand hypergeometric function \cite{AK,G}, was known to, among others, the Japanese authors Mizukawa and Tanaka \cite{MT}, who used them to prove their orthogonality in some special cases.

Later, Mizukawa \cite{M} gave a complete orthogonality proof of \eqref{eq1.19} using Gelfand pairs, then~\cite{MT} used character algebras and closely following this proof came Iliev and Terwilliger's proofs, f\/irst for $n = 2$ \cite{IT}, then for general $n$ in \cite{I}, in which they used tools from Lie algebra theory.  In these proofs the authors found it convenient to use $1 - u_{ij}$ instead of~$u_{ij}$ as parameters in~\eqref{eq1.19}, and also to use
\begin{gather}
\label{eq1.20}
\eta_0 = 1 - \sum_{i=1}^n \eta_i,\qquad m_0 = N - \sum_{i=1}^n m_i,\qquad x_0 = N - \sum_{i=1}^n x_i.
\end{gather}

Our second objective in this paper is to give an alternate proof of orthogonality by using the more elementary method of hypergeometric functions and their various transformation properties and integral representation, and in doing so we f\/ind no special advantage of using $1 - u_{ij}$ instead of $u_{ij}$, or of using \eqref{eq1.20}.  Also, we believe that elementary and perhaps a bit cumbersome as it may be, our method yields a byproduct that throws some light on the underlying geometrical structure of these polynomials, which the other authors may have overlooked.

In Section~\ref{sec2} we will list the transformation and integral representation formulas for the multivariable hypergeometric functions that we shall use throughout this paper.  Section~\ref{sec3} will be devoted to obtaining the necessary conditions of orthogonality, while in Sections~\ref{sec4} and~\ref{sec5} we shall deal with the suf\/f\/icient conditions that will simultaneously establish the orthogonality relation
\begin{gather}
\label{eq1.21}
\underset{\{\bx\}}{\sum} b_n(\bx;N;\bs{\eta})F_1^{(n)}(-\bm;-\bx;-N;\bu) F_1^{(n)}(-\bm';-\bx;-N;\bu) = \d_{\bm,\bm'}/b_n(\bm;N;\overline{\bs{\eta}}),
\end{gather}
where $\overline{\bs{\eta}}$ are the parameters for the dual orthogonality.  The relationship between the $u$'s and the $\eta$'s and $\overline{\eta}$'s will be given in latter sections, specially Section~\ref{sec6}.

In Section~\ref{sec6} we shall examine the geometrical implications of the relationships between the $u$'s that result from the necessary and suf\/f\/icient conditions, while the last section will be aimed at proving that the $n$-dimensional extension of the function in~\eqref{eq1.16} are precisely the eigenfunctions of the kernel $K(\bx;\boy)$. In Section~\ref{sec7} we get the expressions for the
eigenvalues and eigenfunctions of the transition probability kernel governing the evolution of our Markov chain.

Before closing this section we must mention that the polynomials in \eqref{eq1.19} were, in fact, introduced into the statistical literature by R.C.~Grif\/f\/iths~\cite{Gri} as early as 1971 which he def\/ined as coef\/f\/icients in an expansion of their generating function.  However, Mizukawa and Tanaka~\cite{MT} seem to have been the f\/irst to give the explicit expression in~\eqref{eq1.19}.

Few people would disagree with the importance of looking at certain mathematical objects from dif\/ferent
points of view. We feel that this is certainly valid in the case of the present problem: character algebras,
Gelfand--Aomoto functions, Lie algebras and the much older methods of hypergeometric functions including their
series as well as their integral representations have a~useful role to play. Such a wealth of approaches may be
important if one tries to obtain matrix valued versions of these probabilistic models in the spirit of~\cite{Gr,GPT}.
For a very rich and recent extension of the scalar valued solution of the hypergeometric equation, see~\cite{T}. It is
worth noticing that in this case there is yet no extension of the Euler integral representation formula for the $_2F_1$ function, and that algebraic methods such as those coming from Lie algebras or group representation theory, which have
played such in important role in \cite{IT,I,M,M1,MT} have not made a mark in this matrix valued extension yet.
Similar consideration would be relevant if one were to consider a matrix valued extension of, for instance, the work in~\cite{GI}.

\section{Transformation formulas and integral representations}
\label{sec2}

\subsection*{(a) Transformation formulas}

The point of this section is to establish \eqref{eq2.3}, \eqref{eq2.6} and \eqref{eq2.9}.

For references to some of the classical identities in this section the reader can consult for instance \cite{AAR}.

For $n = 0,1,2,\dots$,
\begin{gather*}
{}_2F_1(-n,a;c;x) = \frac {(c-a)_n}{(c)_n} \, {}_2F_1(-n,a;1+a-c-n;1-x).
\end{gather*}

By one iteration,
\begin{gather*}
F_1(-n;a,b;c;x,y) = \frac {(c-a-b)_n}{(c)_n} \, F_1(-n;a,b;1+a+b-c-n;1-x,1-y).
\end{gather*}

By multiple iteration
\begin{gather}
F_1^{(n)}(-m_1,\dots,-m_n;-x_1,\dots,-x_n;-N;u_{11},\dots,u_{1n},u_{21},\dots,u_{2n},u_{n1},\dots,u_{nn}) \nonumber\\
\qquad{} = \frac {\left(\overset{n}{\underset{i=1}{\sum}} x_i-N\right)_{\overset{n}{\underset{i=1}{\sum}} m_i}}{(-N)_{\overset{n}{\underset{i=1}{\sum}} m_i}} F_1^{(n)}\bigg({-}m_1,\dots,-m_n;-x_1,\dots,-x_n; \nonumber\\
\qquad\quad{} N+1-\overset{n}{\underset{i=1}{\sum}} (x_i+m_i);1-u_{11},\dots,1-u_{nn}\bigg) \nonumber\\
\qquad{} = \frac {\left(\overset{n}{\underset{i=1}{\prod}} m_i-N\right)_{\overset{n}{\underset{i=1}{\sum}} x_i}}{(-N)_{\overset{n}{\underset{i=1}{\sum}} x_i}} F_1^{(n)}(\cdots),\label{eq2.3}
\end{gather}
the second identity follows provided the $x$'s are nonnegative integers with $\sum x_i \le N$.

The second transformation for the ${}_2F_1$ function is:
\begin{gather*}
{}_2F_1(a,b;c;x) = (1 - x)^{-a}\, {}_2F_1 \left( a,c-b;c;\frac {-x}{1-x}\right),
\end{gather*}
whose f\/irst iteration gives
\begin{gather*}
F_1(a;b,b';c;x,y)
= (1-y)^{-a}F_1\left(a;b,c-b-b';\frac {x-y}{1-y},\frac {-y}{1-y}\right) \nonumber\\
\phantom{F_1(a;b,b';c;x,y)}{} = (1-x)^{-a}F_1\left(a;c-b-b',b';\frac {-x}{1-x},\frac {y-x}{1-x}\right)
\end{gather*}
and the $n$-th iteration gives
\begin{gather}
F_1^{(n)}(\a_1,\dots,\a_n;\b_1,\dots,\b_n;\g;u_{11},\dots,u_{1n},\dots,u_n,\dots,u_{nn}) \nonumber\\
\qquad{}= (1 - u_{1n})^{-\a_1}(1-u_{2n})^{-\a_2} \cdots (1-u_{nn})^{-\a_n} \nonumber\\
\qquad\quad {}\times  F_1^{(n)}\bigg(\a_1,\dots,\a_n;\b_1,\dots,\b_{n-1},\g-\b_1-\dots -\b_n;\g; \nonumber\\
\qquad \qquad {} \frac {u_{11}-u_{1n}}{1-u_{1n}}, \frac {u_{12}-u_{1n}}{1-u_{1n}}, \dots, \frac {-u_{1n}}{1-u_{1n}},\dots,\frac {u_{1n}-u_{nn}}{1-u_{nn}},\dots,\frac {-u_{nn}}{1-u_{nn}}\bigg),\label{eq2.6}
\end{gather}
which is valid for all $\a$'s, $\b$'s and $\g$, provided the series remain convergent.

\subsection*{(b) Integral representations}

For $0 < \Re b < \Re c$,
\begin{gather*}
{}_2F_1(a,b;c;x) = \frac {\G(c)}{\G(b)\G(c-b)} \int_0^1 \xi^{b-1}(1-\xi)^{c-b-1}(1-\xi x)^{-a}d\xi,
\\
F_1(a;b,c;d;x,y) \nonumber\\
= \frac {\G(d)}{\G(b)\G(c)\G(d - b - c)} \underset{0 < \xi_1+\xi_2 < 1}{\int_0^1\int_0^1}\xi_1^{b-1}\xi_2^{c-1}(1-\xi_1-\xi_2)^{d-b-c-1}(1-\xi_1x-\xi_2y)^{-a}d\xi_1d\xi_2, 
\end{gather*}
provided $\Re(b,c,d,d-b-c) > 0$.

In $n$ dimensions this extends to
\begin{gather}
F_1^{(n)}(\a_1,\dots,\a_n;\b_1,\dots,\b_n;\g;u_{11},\dots,u_{1n},\dots,u_{n1},\dots,u_{nn}) \nonumber\\
\qquad{} = \frac {\G(\g)}{\G\Big(\g-\overset{n}{\underset{i=1}{\sum}}\a_i\Big) \overset{n}{\underset{i=1}{\prod}} \G(\a_i)}
  \underset{0 < \sum\limits_{i=1}^n \xi_i < 1}{\int_0^1 \cdots \int_0^1} \xi_1^{\a_1-1} \cdots \xi_n^{\a_n-1} \left(1 - \overset{n}{\underset{i=1}{\sum}} \xi_i\right)^{\g - \overset{n}{\underset{i=1}{\sum}} \a_i-1} \nonumber\\
\qquad\quad{}\times \overset{n}{\underset{j=1}{\prod}} \left(1 - \overset{n}{\underset{k=1}{\sum}} \xi_ku_{kj}\right)^{-\b_j} d\xi_1\cdots d\xi_n,\label{eq2.9}
\end{gather}
provided
\[
0 < \Re\left(\a_1,\dots,\a_n,\g-\overset{n}{\underset{i=1}{\sum}} \a_i\right).
\]
This is the formula that we shall f\/ind most useful throughout the paper.

\section{Necessary conditions of orthogonality}
\label{sec3}

The point of this section is to show that \eqref{eq3.6} is a necessary condition to insure the desired orthogonality \eqref{eq3.2}.

Denoting the $F_1^{(n)}$ polynomials in \eqref{eq1.19} by $P_{\bm}(\bx)$ for abbreviation, we f\/ind that
\begin{gather}
\label{eq3.1}
\sum_{\{\bx\}} b_n(\bx;N;\bs{\eta})P_{\bm}(\bx) = \prod_{i=1}^n \left( 1 - \sum_{j=1}^n \eta_ju_{ij}\right)^{m_i},
\end{gather}
which is just the generating function of $P_{\bm}(\bx)$.  In a sense this represents the opposite point of view of Grif\/f\/iths~\cite{Gri} where he {\em defined} the polynomials as the coef\/f\/icients of the generating function.

Since one of our aims is to prove the orthogonality relation
\begin{gather}
\label{eq3.2}
I_{\bm}^{\bm'} := \sum_{\{\bx\}} b_n(\bx;N;\bs{\eta})P_{\bm}(\bx)P_{\bm'}(\bx) = 0\qquad \text{if}\quad \bm \ne \bm',
\end{gather}
it must follow, as a necessary condition, that
\begin{gather*}
\sum_{\{\bx\}} b_n(\bx;N;\bs{\eta})P_{\bm}(\bx) = 0,\qquad \bm \ne (0,\dots,0),
\end{gather*}
which, by \eqref{eq3.1}, implies that
\[
\prod_{i=1}^n \left(1 - \sum_{j=1}^n \eta_ju_{ij}\right)^{m_i} = 0,
\]
and therefore
\begin{gather}
\label{eq3.4}
\sum_{j=1}^n \eta_ju_{ij} = 1,\qquad i = 1,2,\dots,n.
\end{gather}
If we denote the parameters of the dual orthogonality by $\overline{\bs{\eta}} = (\overline{\eta}_1,\overline{\eta}_2,\dots,\overline{\eta}_n)$, then, using
\begin{gather*}
\sum_{\{\bm\}} b_n(\bm;N;\overline{\bs{\eta}})P_{\bm}(\bx) = 0,\qquad \bx \ne (0,0,\dots,0),
\end{gather*}
we get as a necessary condition
\begin{gather}
\label{eq3.6}
\sum_{j=1}^n \overline{\eta}_ju_{ji} = 1,\qquad i = 1,\dots,n.
\end{gather}

\section{Suf\/f\/icient conditions of orthogonality}
\label{sec4}

In this section we show that the relation \eqref{eq4.11} among the parameters is
suf\/f\/icient for orthogonality.

Instead of the sum in \eqref{eq3.2} let us consider, for the time being, the following sum
\begin{gather*}
\!\sum_{\{\bx\}} b_n(\bx;N;\bs{\eta})P_{\bm}(x)F_1^{(n)}\big(\a_1,\a_2,\dots,\a_n;-m_1,\dots,-m_n;\g;
u_{11},\dots,u_{1n},\dots,u_{n1},\dots,u_{nn}\big),\!
\end{gather*}
in which we assume that while the $m$'s are nonnegative integers the $\a$'s are not, and nor is $\g$ a~nonpositive integer, and that
\begin{gather*}
0 < \Re\left( \a_1,\dots,\a_n,\g - \sum_{i=1}^n \a_i\right).
\end{gather*}

Using \eqref{eq2.9} we then f\/ind that the above sum equals
\begin{gather}
\frac {\G(\g)}{\G\left(\g - \overset{n}{\underset{i=1}{\sum}} \a_i\right)\overset{n}{\underset{i=1}{\prod}} \G(\a_i)} \underset{0 < \overset{n}{\underset{i=1}{\sum}} \xi_i < 1}{\int_0^1 \cdots \int_0^1} \left( 1 - \overset{n}{\underset{i=1}{\sum}} \xi_i\right)^{\g-\overset{n}{\underset{i=1}{\sum}} \a_i - 1} \nonumber\\
\qquad{}\times \prod_{i=1}^n \xi_i^{\a_i-1} d\xi_i \sum_{\{\bx\}} b_n(\bx;N;\bs{\eta})P_{\bm}(\bx)
 \prod_{j=1}^n \left( 1 - \sum_{i=1}^n \xi_iu_{ij}\right)^{x_j}.\label{eq4.3}
\end{gather}
However,
\begin{gather*}
 \sum_{\{\bx\}} b_n(\bx;N;\bs{\eta}) \prod_{j=1}^N \left( 1 - \sum_{i=1}^n \xi_iu_{ij}\right)^{x_j} \prod_{i=1}^n (-x_i)_{\overset{n}{\underset{j=1}{\sum}} r_{ji}} \nonumber\\
\qquad{} = (-N)_{\underset{i,j}{\sum} r_{ij}} \left( 1 - \sum_{i,j} \xi_i\eta_ju_{ij}\right)^{N-\underset{i,j}{\sum} r_{ij}}
  \prod_{k=1}^n \left\{ \eta_k\left( 1 - \sum_{i=1}^n \xi_iu_{ik}\right)\right\}^{\overset{n}{\underset{j=1}{\sum}} r_{jk}}.
\end{gather*}
But, by \eqref{eq3.4}, $\sum\limits_{j=1}^n \eta_ju_{ij} = 1$, so $\sum\limits_{i,j} \xi_i\eta_ju_{ij} = \sum\limits_{i=1}^n \xi_i$, and hence the sum inside the integral in~\eqref{eq4.3} becomes
\begin{gather}
 \left( 1 - \sum \xi_i\right)^N \! \sum_{\{r_{ij}\}} \frac {\overset{n}{\underset{i=1}{\sum}} (-m_i)_{\overset{n}{\underset{j=1}{\sum}} r_{ij}}}{\underset{i,j}{\prod} r_{ij}!} \!\left( 1 - \sum_{i=1}^n \xi_i\right)^{-\underset{i,j}{\sum} r_{ij}}   \prod_{i,j} u_{ij}^{r_{ij}} \prod_{k=1}^n \left\{ \eta_k \left( 1 - \sum_{i=1}^n \xi_iu_{ik}\right)^{\overset{n}{\underset{j=1}{\sum}}r_{jk}} \right\} \nonumber\\
\qquad{}= \left( 1 - \sum_{i=1}^n \xi_i\right)^{N - \overset{n}{\underset{k=1}{\sum}} m_k} \prod_{r=1}^n \left\{ \D_r\xi_r + \sum_{s\ne r}^n A_{r,s}\xi_s\!\right\}^{m_r},\label{eq4.5}
\end{gather}
where
\begin{gather}
\D_r  = \sum_{s=1}^n \eta_su_{rs}^2 - 1,\qquad r = 1,2,\dots,n, \nonumber\\
A_{r,s}  = \sum_{j=1}^n \eta_ju_{rj}u_{sj} - 1  = A_{s,r}   ,\qquad r,s = 1,2,\dots,n.\label{eq4.6}
\end{gather}

We now expand the $n$-fold product on the right-hand side of \eqref{eq4.5} to get
\begin{gather}
 \sum_{\{k_i\}} \prod_{i=1}^n \binom{m_i}{k_{i,1},k_{i,2},\dots,k_{i,n-1}}
  \prod_{i=1}^N \D_i^{m_i - \overset{n-1}{\underset{j=1}{\sum}} k_{ij}}
  A_{1,2}^{k_{11}+k_{21}} A_{1,3}^{k_{12}+k_{21}} \cdots A_{1,n}^{k_{1,n-1} + k_{n,n-1}} \nonumber\\
\qquad{}\times A_{2,3}^{k_{22}+k_{32}} A_{2,4}^{k_{23}+k_{42}} \cdots A_{n,n-1}^{k_{n,n-1}+k_{n-1,n-1}}
  \xi_1^{m_1-\overset{n-1}{\underset{j=1}{\sum}} k_{1,j}+(k_{2,1}+k_{n,1})} \nonumber\\
\qquad{}\times \xi_2^{m_2 - \overset{n-1}{\underset{j=1}{\sum}} k_{2,j}+(k_{11}+k_{32}+\cdots+k_{n2})}
\cdots  \xi_n^{m_n - \overset{n-1}{\underset{j=1}{\sum}} k_{n,j} + (k_{1,n-1}+k_{2,n-1}+\cdots+k_{n-1,n-1})}.\label{eq4.7}
\end{gather}

Substitution of \eqref{eq4.5} and \eqref{eq4.7} inside the integrand of \eqref{eq4.3} and computing the integral, with the $\a_i$'s replaced by $-m_i$'s and $\g$ by $-N$, gives the value of the integral as
\begin{gather}
 (-m'_1)_{m_1 + \overset{n}{\underset{r=2}{\sum}} k_{r1} - \overset{n-1}{\underset{j=1}{\sum}} k_{1j}}(-m'_2)_{m_2+k_{11}+\overset{n}{\underset{r=3}{\sum}} k_{r2} - \overset{n-1}{\underset{j=1}{\sum}} k_{2j}} \nonumber\\
\qquad{}\times \cdots \times (-m'_n)_{m_n + \overset{n-1}{\underset{r=1}{\sum}} k_{r,n-1} - \overset{n-1}{\underset{j=1}{\sum}} k_{n,j}} (m'_1+\cdots+m'_n-N)_{N-m_1-\cdots-m_n} \nonumber\\
\qquad{}\times  A_{12}^{k_{11}+k_{21}} \cdots A_{n,n-1}^{k_{n,n-1}+k_{n-1,n-1}}   \prod_{r=1}^n \D_r^{m_r - \overset{n-1}{\underset{j=1}{\sum}} k_{rj}}\Big/ (-N)_N.\label{eq4.8}
\end{gather}
Let $\bm = \boe_i$, $\bm' = \boe_j$, where $\boe_i = (0,\dots,1,0,\dots,0)$, $\boe_j = (0,\dots,0,\dots,1,\dots,0)$, i.e., $e$'s are unit vectors with the $i$-th component being $1$ for $\boe_i$ and $j$-th position for $\boe_j$, otherwise all components~$0$.

In particular, for $\bm = \boe_1$, and $\bm' = \boe_n$,
\begin{gather*}
 -\frac {1}{N} \sum_{\{k_{1,j}\}} \binom{1}{k_{11},\dots,k_{1,n-1}} \D_1^{1-\overset{n-1}{\underset{i=1}{\sum}} k_{ij}} A_{12}^{k_{11}}A_{13}^{k_{12}} \cdots A_{1n}^{k_{1,n-1}} \nonumber\\
 \qquad{}\times (0)_{1-\overset{n-1}{\underset{j=1}{\sum}} k_{1j}}(0)_{k_{11}}(0)_{k_{1,2}}\cdots (-1)_{k_{1,n-1}}
 = -\frac {A_{1,n}}{N} .
\end{gather*}

By similar arguments it follows that
\begin{gather*}
A_{i,j} = 0,\qquad i \ne j,
\end{gather*}
which, by \eqref{eq4.6}, means that
\begin{gather}
\label{eq4.11}
\sum_{j=1}^n \eta_ju_{rj}u_{sj} = 1,\qquad r,s = 1,2,\dots,n,\qquad r \ne s.
\end{gather}

So the only terms that survive in \eqref{eq4.8} give
\begin{gather}
\label{eq4.12}
I_{\bm}^{\bm'} = \frac {\left(-N + \overset{n}{\underset{k=1}{\sum}} m'_k\right)_{N - \overset{n}{\underset{j=1}{\sum}} m_j} \overset{n}{\underset{k=1}{\prod}} (-m'_k)_{m_k}}{(-N)_N} \prod_{j=1}^n \D_j^{m_j}.
\end{gather}

The f\/irst factor in the numerator on the right-hand side of \eqref{eq4.12} is zero unless $N - \overset{n}{\underset{k=1}{\sum}} m'_k \ge N - \overset{n}{\underset{k=1}{\sum}} m_k$, i.e.,
\begin{gather}
\label{eq4.13}
\sum_{j=1}^n m_j \ge \sum_{j=1}^n m'_j,
\end{gather}
while the $n$ remaining factors imply that they are zero unless
\begin{gather*}
m'_i \ge m_i,\qquad i = 1,2,\dots,n.
\end{gather*}
\eqref{eq4.12} and \eqref{eq4.13} imply that $I_{\bm}^{\bm'} = 0$ unless $\bm' = \bm$.  So we have
\begin{gather*}
I_{\bm}^{\bm'} = \frac {\overset{n}{\underset{k=1}{\sum}} (-m_k)_{m_k}}{(-N)_{\overset{n}{\underset{k=1}{\sum}} m_k}} \prod_{j=1}^n \D_j^{m_j} \d_{\bm,\bm'},
\end{gather*}
and that \eqref{eq4.11} is a suf\/f\/icient condition for orthogonality.  Clearly, it can be written in a form closer to the desired form \eqref{eq1.21}, namely
\begin{gather}
\label{eq4.16}
I_{\bm}^{\bm'} = \frac {\overset{n}{\underset{j=1}{\prod}} \D_j^{m_j}}{\binom{N}{m_1,m_2,\dots,m_n}} \d_{\bm,\bm'}.
\end{gather}

In the following section we complete this part of the work by expressing $\overset{n}{\underset{j=1}{\prod}} \D_j^{m_j}$ in terms of the parameters of dual orthogonality, i.e., $\overline{\eta}_i$'s.

\section{Reduction of (\ref{eq4.16}) to (\ref{eq1.21})}
\label{sec5}

As a typical $\D$ let us consider
\begin{gather*}
\D_1 = \sum_{s=1}^n \eta_su_{1s}^2 - 1 \qquad \text{by \eqref{eq4.6}.}
\end{gather*}
From \eqref{eq3.4} we have
\begin{gather}
\label{eq5.2}
\sum_{j=1}^n \eta_ju_{1j} = 1,
\end{gather}
while \eqref{eq4.11} gives
\begin{gather}
\label{eq5.3}
\sum_{j=1}^n \eta_ju_{1j}u_{rj} = 1,\qquad r = 2,\dots,n.
\end{gather}
From \eqref{eq5.2} and \eqref{eq5.3} it follows that
\begin{gather*}
\eta_1  = \frac {u_{12}u_{13}\cdots u_{1n}}{\L}
\begin{vmatrix}
1 & 1 & \dots & 1 \\
1 & u_{22} & \dots & u_{2n} \\
\vdots & \vdots \\
1 & u_{n2} & \dots & u_{nn}
\end{vmatrix} , \nonumber\\
\eta_r  = \frac {u_{11}\dots u_{1,r-1} u_{1,r+1}\cdots u_{1n}}{\L}
\begin{vmatrix}
1 & \dots & 1 & 1 & \dots & 1 \\
u_{11} & \dots & u_{2,r-1} & 1 & \dots & u_{2n} \\
\vdots \\
k_{n1} & \dots & u_{n,r-1} & 1 & \dots & u_{nn}
\end{vmatrix} ,\qquad r = 2,\dots,n,
\end{gather*}
where
\begin{gather*}
\L = \prod_{k=1}^n u_{1k} \begin{vmatrix}
1 & 1 & \dots & 1 \\
u_{21} & u_{22} & \dots & u_{2n} \\
\vdots \\
u_{n1} & u_{n2} & \dots & u_{nn}
\end{vmatrix} = D\overline{\eta}_1 \prod_{k=1}^n u_{1k},
\qquad
D = \Det(u_{ij}),
\end{gather*}
which follows from \eqref{eq3.4}.  So
\begin{gather}
\left( \sum_{s=1}^n \eta_su_{1s}^2 - 1\right) D\overline{\eta}_1  = u_{11} \begin{vmatrix}
1 & 1 & \dots & 1 \\
1 & u_{22} & & u_{2n} \\
\vdots & \vdots & & \vdots \\
1 & u_{n2} & & u_{nn}
\end{vmatrix} + u_{12} \begin{vmatrix}
1 & 1 & \dots & 1 \\
u_{21} & 1 & \dots & u_{2n} \\
\vdots \\
u_{n1} & 1 & \dots & u_{nn}
\end{vmatrix}  + \cdots \nonumber\\
 \hphantom{\left( \sum_{s=1}^n \eta_su_{1s}^2 - 1\right) D\overline{\eta}_1  =}{}
 + u_{1n} \begin{vmatrix}
1 & \dots & 1 & 1 \\
u_{21} & \dots & u_{2,n-1} & 1 \\
\vdots \\
u_{n1} & \dots & u_{n,n-1} & 1
\end{vmatrix} - D\overline{\eta}_1.\label{eq5.7}
\end{gather}

It is clear that expansion by the f\/irst rows of the $n$ determinants on the right produces $n$ terms:
\[
D + u_{11} \begin{vmatrix}
0 & 1 & \dots & 1 \\
1 & u_{22} & \dots & k_{2n} \\
\vdots \\
1 & u_{n2} & \dots & u_{nn}
\end{vmatrix} + u_{12} \begin{vmatrix}
1 & 0 & 1 & \dots & 1 \\
u_{21} & 1 & u_{23} & \dots & u_{2n} \\
\vdots \\
u_{n1} & 1 & u_{n3} & \dots & u_{nn}
\end{vmatrix} + \dots + u_{1n} \begin{vmatrix}
1 & \dots & 1 & 0 \\
u_{21} & \dots & u_{n,n-1} & 1 \\
\vdots \\
u_{n1} & \dots & u_{n,n-1} & 1
\end{vmatrix}
\]
whence, we derive by careful rearrangement of the determinants the following expression:
\begin{gather}
D - \begin{vmatrix}
u_{11} & u_{12} & \dots & u_{1n} \\
1 & 1 & \dots & 1 \\
u_{31} & u_{32} & \dots & u_{3n} \\
\vdots \\
u_{n1} & u_{n2} & \dots & u_{nn}
\end{vmatrix} - \dots - \begin{vmatrix}
u_{11} & u_{12} & \dots & u_{1n} \\
\vdots \\
\vdots \\
u_{n-1,1} & u_{n-1,2} & \dots & u_{n-1,n} \\
1 & 1 & \dots & 1
\end{vmatrix} \nonumber\\
\qquad{} = D(1 - \overline{\eta}_2 - \cdots - \overline{\eta}_n).\label{eq5.8}
\end{gather}
From \eqref{eq5.7} and \eqref{eq5.8} we f\/inally obtain the desired expression
\begin{gather*}
\sum_{s=1}^n \eta_su_{1s}^2 - 1 = \frac {1 - \overset{n}{\underset{s=1}{\sum}} \overline{\eta}_s^2}{\overline{\eta}_1} ,
\end{gather*}
and, in a similar way, we derive the general form:
\begin{gather*}
\sum_{s=1}^n \eta_su_{rs}^2 - 1 = \frac {1 - \overset{n}{\underset{s=1}{\sum}} \overline{\eta}_s^2}{\overline{\eta}_r} ,
\qquad
r = 1,2,\dots,n.
\end{gather*}

\section{Relations among the parameters}
\label{sec6}

The previous sections were devoted to proving the orthogonality of the $n^2$-parameter polynomials $P_{\bm}(x)$ def\/ined in \eqref{eq1.19}, with respect to the multinomial distribution $b_n(\bx;N;\bs{\eta})$, $\bs{\eta} = (\eta_1,\dots,\eta_n)$ (or $b_n(\bm;N;\overline{\bs{\eta}})$, $\overline{\bs{\eta}} = (\overline{\eta}_1,\dots,\overline{\eta}_n)$, with $\sum \eta_i = \sum \overline{\eta}_i$).  Clearly, the $n$ relations in \eqref{eq3.4} completely def\/ine the $\eta_i$'s in terms of the $u$'s, as they do the $\overline{\eta}_i$'s, by \eqref{eq3.6}.  If we think of the $\eta$'s as given parameters then the $n^2$ polynomial parameters must satisfy in addition $n^2 - (2n-1) = (n-1)^2$ relations among themselves.  Where do they come from and what do they mean geometrically is the question we shall examine in this section.

Note that the suf\/f\/icient conditions of orthogonality \eqref{eq4.11} give $\binom{n}{2}$ further relations between the $\eta$'s and $u$'s.  In fact, combining \eqref{eq3.4} and \eqref{eq4.11} we get
\begin{gather}
\label{eq6.1}
\sum_{j=1}^n \eta_j u_{rj}(1 - u_{sj}) = 0,\qquad r,s = 1,2,\dots,n,\quad r \ne s.
\end{gather}
Clearly, the positivity of all the probability parameters, i.e., $\eta$'s, require that the determinant $\Det(u_{rj}(1 - u_{sj})) = 0$.  It is possible to analyze these determinants and obtain the missing $(n-1)^2$ relations.  We shall, however, take a dif\/ferent route.

Our proof of orthogonality is based almost entirely on the integral representation \eqref{eq2.9} and the simple multinomial summation formula.  The proof is direct and elementary.  However, if we had instead used the transformations \eqref{eq2.3} and \eqref{eq2.6} we could have reduced the problem to a~$(n-1)$-variable case, with the polynomials having $(n-1)^2$ parameters, instead of $n^2$.  There are, of course, $n$ ways to make this reduction, depending on which of the $n$ variables we integrate f\/irst.  If we do the $x_n$ summation f\/irst, then by $n$ long and tedious set of computations we can f\/ind that the $I_{\bm}^{\bm'}$ reduces to an expression proportional to the $(n-1)$-fold sum
\begin{gather*}
 \sum_{\{x_1,\dots,x_{n-1}\}} \binom{\overset{n}{\underset{i=1}{\sum}} m_i}{x_1,\dots,x_{n-1}} \zeta_1^{x_1} \zeta_2^{x_2} \cdots \zeta_{n-1}^{x_{n-1}} \nonumber\\
{}\times F^{(n-1)} \Bigg( {-}m'_1,\dots,-m'_{n-1};-x_1,\dots,-x_{n-1};-\overset{n}{\underset{i=1}{\sum}} m_i;1 - \frac {u_{nn}u_{11}}{u_{n1}u_{1n}},1 - \frac {u_{nn}u_{12}}{u_{1n}u_{n2}},\dots,  \nonumber \\
 1 - \frac {u_{nn}u_{1,n-1}}{u_{1n}u_{n,n-1}},1 - \frac {u_{nn}u_{21}}{u_{1n}u_{n2}},\dots,1 - \frac {u_{nn}u_{2,n-1}}{u_{1n}u_{n,n-1}},\dots,
  1 -   \frac {u_{nn}u_{n-1,1}}{u_{1n}u_{n,n-1}},\dots,1 - \frac {u_{nn}u_{n-1,n-1}}{u_{1n}u_{n,n-1}}\Bigg) \nonumber\\
{}\times F^{(n-1)} \Bigg( {-}m_1,\dots,-m_{n-1};-x_1,\dots,-x_{n-1};-\overset{n}{\underset{i=1}{\sum}} m_i;1 - \frac {(1-u_{nn})(1-u_{11})}{(1-u_{n1})(1-u_{1n})}, \nonumber\\
1 -   \frac {(1-u_{nn})(1-u_{1n})}{(1-u_{1n})(1-u_{n2})},\dots,1 - \frac {(1-u_{nn})(1-u_{1,n-1})}{(1-u_{1n})(1-u_{n,n-1})},1 - \frac {(1-u_{nn})(1-u_{21})}{(1-u_{2n})(1-u_{n2})},\dots\Bigg), 
\end{gather*}
where
\begin{gather*}
\zeta_1  = \frac {\eta_1u_{n1}(1-u_{n1})}{\eta_nu_{nn}(1-u_{nn})},\quad  \zeta_2 = \frac {\eta_2u_{n2}(1-u_{n2})}{\eta_nu_{nn}(1-u_{nn})} ,\quad \dots, \quad
\zeta_n  = \frac {\eta_{n-1}u_{n,n-1}(1-u_{n,n-1})}{\eta_nu_{nn}(1-u_{nn})} . 
\end{gather*}

If we were to have orthogonality of the $n$-variable case then surely a reduction to a lower dimension will retain the same property.  However, the appearance of the two $F_1^{(n-1)}$ functions above doesn't suggest their parameters are the same.  But the point is that they are, not identically, but consistent with the necessary conditions of orthogonality~\eqref{eq3.4}.  To illustrate this point let us take the product of $\zeta$'s and the f\/irst $n-1$ parameters of the f\/irst $F_1^{(n-1)}$, i.e., compute
\begin{gather*}
 \zeta_1\left(1 - \frac {u_{nn}u_{11}}{u_{n1}-u_{1n}}\right) + \zeta_2 \left(1 - \frac {u_{nn}u_{12}}{u_{1n}u_{n2}}\right) + \cdots
 + \zeta_n\left(1 - \frac {u_{nn}u_{1,n-1}}{u_{1n}u_{n,n-1}}\right)  \nonumber\\
 \qquad{} = \sum_{i=1}^{n-1} \zeta_i - \frac {\eta_1u_{11}(1-u_{n1})}{\eta_nu_{1n}(1-u_{nn})} - \frac {\eta_2u_{12}(1-u_{n2})}{\eta_nu_{1n}(1-u_{nn})} - \cdots - \frac {\eta_{n-1}u_{1,n-1}(1-u_{n,n-1})}{\eta_nu_{1n}(1 - u_{nn})}
\\ 
 \qquad{}
 = \sum_{i=1}^{n-1} \zeta_i - \frac {\overset{n-1}{\underset{j=1}{\sum}} \eta_ju_{1j}(1 - u_{nj})}{\eta_nu_{1n}(1-u_{nn})}
 = \sum_{i=1}^{n-1} \zeta_i + \frac {\eta_nu_{1n}(1-u_{nn})}{\eta_nu_{nn}(1-u_{nn})} = \sum_{i=1}^{n-1} \zeta_i+1 \qquad \text{by \eqref{eq6.1}.}\nonumber
\end{gather*}
On the other hand,
\begin{gather*}
 \zeta_1 \left(1 - \frac {(1-u_{nn})(1-u_{11})}{(1-u_{n1})(1-u_{1n})}\right) + \zeta_2 (1 - \frac {(1-u_{nn})(1-u_{12})}{(1-u_{1n})(1-u_{n2})}) + \cdots
  \nonumber\\
  \qquad{} + \zeta_{n-1} \left(1 - \frac {(1-u_{nn})(1-u_{n,n_2})}{(1-u_{1n})(1-u_{n,n2})}\right)
 = \sum_{i=1}^{n-1} \zeta_i - \frac {\eta_1u_{n1}(1-u_{11})}{\eta_nu_{nn}(1-u_{1n})} \nonumber\\
 \qquad{} - \frac {\eta_2u_{n2}(1-u_{12})}{\eta_nu_{nn}(1-u_{1n})} - \cdots - \frac {\eta_{n-1}u_{n,n-1}(1-u_{1,n-1})}{\eta_nu_{nn}(1-u_{1,n})}
 = \sum_{i=1}^{n-1} \zeta_i + 1.
\end{gather*}
Similarly the equalities of other parameters are also established.  The $(n-1)^2$ relations can be expressed in the following compact form:
\begin{gather}
\label{eq6.6}
U_{jk}U_{nn} = U_{nk}U_{jn},
\end{gather}
where
\begin{gather*}
U_{ij} = 1 - u_{ij}^{-1},
\end{gather*}
(one must, of course, take it for granted that $u_{ij} \ne 0$, for any $i$, $j$).

The geometrical implication of these relations seems to suggest that the $n$-variable Kraw\-tchouk polynomials \eqref{eq1.19} live on an $(n-1)^2$-dimensional submanifold, def\/ined by \eqref{eq6.6}, of the $n^2$-dimensional space of $(u_{ij})$.

\section[Eigenvalues and eigenfunctions of $K(\bi;\bj)$]{Eigenvalues and eigenfunctions of $\boldsymbol{K(\bi;\bj)}$}
\label{sec7}

For the sake of notational consistency and convenience we shall replace $K(\bi;\bj)$ by $K(\bx;\boy)$, so that by use of \eqref{eq1.10} and \eqref{eq1.15}, our eigenvalue problem can be stated as follows:  f\/ind $\psi_{\bm}(\bx)$ such that
\begin{gather}
 b_n(\bx;N;\bs{\beta}) \sum_{\{\boy\}} b_n(\boy;N;\bs{\eta}) \prod_{k=1}^n (1-\a_k)^{y_k}\sum_{\{\br\}} \frac {(-x_1)_{r_1}\cdots(-x_n)_{r_n}(-y_1)_{r_1}\cdots(-y_n)_{r_n}}{r_1!\cdots r_n!(-N)_{r_1+r_2+\cdots+r_n}} \nonumber \\
\qquad{}\times \left( \frac {\a_1}{\b_1(\a_1-1)}\right)^{r_1} \cdots
 \left(\frac {\a_n}{\b_n(\a_n-1)}\right)^{r_n} \psi_{\bm}(\boy) = \l_{\bm}b_n(\bx;N;\bs{\eta})\psi_{\bm}(\bx).\label{eq7.1}
\end{gather}
We will show that
\begin{gather*}
\psi_{\bm}(\bx) = P_{\bm}(\bx)
\end{gather*}
for some choices of the parameters $u_{ij}$'s.  To this end what we will do is compute the sum above with $\psi_{\bm}(\bx)$ replaced by $P_{\bm}(\bx)$.

Using the integral representation \eqref{eq2.9} we f\/ind that the sum that needs to be computed is
\begin{gather}
 \sum_{\{\boy\}} \bb(\boy;N;\bs{\eta}) \prod_{k=1}^n (-y_k)_{r_k} \prod_{j=1}^n \left\{ (1 - \a_j)\left( 1 - \sum_{i=1}^n \xi_iu_{ij}\right)\right\}^{y_k}
 \nonumber\\
 \qquad{} = (-N)_{\overset{n}{\underset{i=1}{\sum}} r_i} \prod_{j=1}^n \left\{ \eta_j(1-\a_j)\left( 1 - \sum_{i=1}^n \xi_iu_{ij}\right)\right\}^{r_j}
 \nonumber\\
\qquad\quad{}\times \left[ 1 - \sum_{i=1}^n \eta_i + \sum_{j=1}^n \eta_i(1-\a_j)\left(1 - \sum_{i=1}^n \xi_iu_{ij}\right)\right]^{N - \overset{n}{\underset{i=1}{\sum}} r_i}.\label{eq7.3}
\end{gather}
Since by \eqref{eq1.13}
\begin{gather*}
\sum_{j=1}^n \eta_j(1-\a_j) = \frac {1 - \sum\limits_{i=1}^n \eta_i}{1 - \sum\limits_{i=1}^n \b_i} \sum_{j=1}^n \b_j,
\end{gather*}
we have
\begin{gather*}
1 - \sum_{i=1}^n \eta_i + \sum_{j=1}^n \eta_j (1-\a_j) = \frac {1 - \sum\limits_{i=1}^n \eta_i}{1 - \sum\limits_{i=1}^n \b_i} .
\end{gather*}

Therefore the right-hand side of \eqref{eq7.3} simplif\/ies to
\begin{gather*}
 (-N)_{r_1+\cdots+r_n} \left( \frac {1 - \overset{n}{\underset{i=1}{\sum}} \eta_i}{1 - \overset{n}{\underset{i=1}{\sum}} \b_i} \right)^{N-\overset{n}{\underset{i=1}{\sum}} r_i} \left[ 1 - \frac {1 - \sum\limits_{i=1}{n} \b_i}{1 - \sum\limits_{i=1}^n \eta_i} \sum_i \sum_j \xi_iu_{ij}\eta_j(1 - \a_j)\right]^{N - \overset{n}{\underset{i=1}{\sum}} r_i} \nonumber\\
\qquad\quad{}  \times \prod_{j=1}^n \left\{ \eta_j(1-\a_j)\left( 1 - \sum_{i=1}^n \xi_iu_{ij}\right)\right\}^{r_j} \nonumber\\
 \qquad{} = (-N)_{\sum\limits_{i=1}^n r_i} \left( \frac {1 -  \sum\limits_{i=1}^n \eta_i}{1 -  \sum\limits_{i=1}^n \b_i} \right)^N \left[ 1 - \frac {1 -  \sum\limits_{i=1}^n \b_i}{1 -  \sum\limits_{i=1}^n \eta_i} \sum_i v_i\xi_i\right]^{N -  \sum\limits_{i=1}^n r_i}\nonumber\\
\qquad\quad{}  \times
  \prod_{j=1}^n \left\{ \b_j \left( 1 -  \sum\limits_{i=1}^n \xi_iu_{ij}\right)\right\}^{r_j},\\ 
v_i =  \sum\limits_j \eta_j(1-\a_j)u_{ij} = \frac {1-\sum\limits_{k=1}^n \eta_k}{1-\sum\limits_{k=1}^n \b_k} \sum\limits_j \b_ju_{ij} .\nonumber
\end{gather*}
So the sum over the $r_i$'s is
\begin{gather*}
\left( \frac {1 -  \sum\limits_{i=1}^n \eta_i}{1 - \sum\limits_{i=1}^n \b_i} \right)^N \sum_{\{r_i\}} \prod_{k=1}^n \binom{x_k}{r_k} \prod_{j=1}^n \left\{ \frac {\a_j}{1-\a_j} \left( 1 - \sum\limits_{i=1}^n\xi_iu_{ij}\right)\right\}^{r_j}
  \left\{ 1 - \sum\limits_{i=1}^n \omega_i\xi_i\right\}^{N - \sum\limits_{i=1}^n r_i} \nonumber\\
  = \left( \frac {1 - \sum\limits_{i=1}^n \eta_i}{1 - \sum\limits_{i=1}^n \b_i}\right)^N \! \left( 1 - \sum\limits_{i=1}^n \omega_i\xi_i\right)^{N - \sum\limits_{i=1}^n x_i} \! \prod_{j=1}^n \left\{ \frac {\a_j}{1-\a_j} + 1 - \sum\limits_{i=1}^n \omega_i\xi_i - \frac {\a_j}{1-\a_j} \sum_{j=1}^n \xi_iu_{ij}\right\}^{x_j} \nonumber\\
= \left( \frac {1 - \sum\limits_{i=1}^n \eta_i}{1 - \sum\limits_{i=1}^n \b_i} \right)^N \prod_{j=1}^n (1-\a_j)^{-x_j}
  \prod_{j=1}^n \left\{ 1 - \a_j \sum_{i=1}^n \xi_iu_{ij} - (1-\a_j) \sum\limits_{i=1}^n \omega_i\xi_i\right\}^{x_j}\nonumber \\
\quad{} \times \left( 1 - \sum_{i=1}^n \omega_i\xi_i\right)^{N - \sum\limits_{i=1}^n x_i}
\end{gather*}
with
\[
\omega_i = \sum_{j=1}^n \b_ju_{ij}.
\]
Since
\begin{gather*}
\prod_{j=1}^n (1-\a_j)^{-x_j} = \left( \frac {1 - \sum\limits_{i=1}^n \b_i}{1 - \sum\limits_{i=1}^n \eta_i}\right)^{\sum\limits_{i=1}^n x_i} \prod_{j=1}^n \left( \frac {\eta_j}{\b_j} \right)^{-x_j},
\end{gather*}
the right-hand side of \eqref{eq7.1}, via \eqref{eq2.9}, becomes
\begin{gather}
 b_n(\bx;N;\bs{\eta}) \underset{0 < \sum\limits_{i=1}^n \xi_i < 1}{\int_0^1 \cdots \int_0^1} \xi_1^{a_i-1} \cdots \xi_n^{a_n-1} \left( 1 - \sum\limits_{i=1}^n  \xi_i\right)^{c - \sum\limits_{i=1}^n a_i - 1} \nonumber\\
\qquad\quad{}\times \prod_{j=1}^n \left\{ 1 - \sum_{i=1}^n (\a_ju_{ij} + (1 - \a_j)\omega_i)\xi_i\right\}^{x_j}
  \left( 1 - \sum_{i=1}^n \omega_i\xi_i\right)^{N - \sum\limits_{i=1}^n x_i} \prod_{j=1}^n d\xi_j \nonumber\\
\qquad {} = b_n(\bx;N;\bs{\eta}) \sum_{\underset{0 \le \underset{i,j}{\sum} k_{ij} \le N}{\{k_{ij}\}}} \frac {\overset{n}{\underset{i=1}{\prod}} (-m_i)_{\overset{n+1}{\underset{j=1}{\sum}} k_{ij}} \overset{n}{\underset{i=1}{\prod}} (-x_i)_{\overset{n}{\underset{j=1}{\sum}} k_{ji}}}{\underset{i,j}{\prod} k_{ij}!(-N)_{\overset{n}{\underset{i=1}{\sum}} \overset{n+1}{\underset{j=1}{\sum}} k_{i,j}}} \nonumber\\
\qquad\quad{}\times (x_1+\cdots+x_n-N)_{\overset{n}{\underset{i=1}{\sum}} k_{i,n+1}}   \prod_{i=1}^n \omega_i^{k_{i,n+1}}
  \prod_{i=1}^n \prod_{j=1}^{n+1} \{\a_ju_{ij} + (1 - \a_j)\omega_i\}^{k_{ij}} \nonumber\\
\qquad {} = b_n(\bx;N;\bs{\eta}) \prod_{i=1}^N (1 - \omega_i)^{m_i}
  \sum_{\underset{0 \le \underset{k_{ij}}{\sum} \le N}{\{k_{ij}\}}} \frac {\overset{n}{\underset{i=1}{\prod}} (-m_i)_{\overset{n+1}{\underset{j=1}{\sum}} k_{ij}} \overset{n}{\underset{i=1}{\prod}} (-x_i)_{\overset{n}{\underset{j=1}{\sum}} k_{ji}}}{\underset{i,j}{\prod} k_{ij}!(-N)_{\overset{n}{\underset{i,j}{\sum}} k_{ij}}} \nonumber\\
\qquad\quad{}\times \prod_{i,j=1}^n \left( \frac {\a_ju_{ij} + (1 - \a_j)\omega_i - \omega_i}{1-\omega_i}\right)^{k_{ij}}.\label{eq7.9}
\end{gather}

So $P_{\bm}(\bx)$ is an eigenfunction of $K(\bx;\boy)$, provided $u_{ij}$ exists for all $i,j = 1,2,\dots,n$, such that
\begin{gather*}
\a_j(u_{ij}-\omega_i) = (1 - \omega_i)u_{ij},
\end{gather*}
i.e.,
\begin{gather}
\a_j(1 - \omega_i) - \a_j(1 - u_{ij}) = (1 - \omega_i)u_{ij},\nonumber\\
\label{eq7.11}
\frac {\a_j(1 - u_{ij})}{\a_j-u_{ij}} = 1 - \omega_i = 1 - \sum_{k\ne j}^n \b_ku_{ik} - \b_ju_{ij}.
\end{gather}
So
\begin{gather*}
\sum_{k\ne j}^n \b_k u_{ik}  = 1 - \b_ju_{ij} - \frac {\a_j(1-u_{ij})}{\a_j-u_{ij}}
 = -\b_ju_{ij} + \frac {u_{ij}(1 - \a_j)}{u_{ij}-\a_j}  \nonumber\\
\hphantom{\sum_{k\ne j}^n \b_k u_{ik}}{} = \frac {\b_ju_{ij}^2 - (1 - \a_j + \a_j\b_j)u_{ij}}{\a_j-u_{ij}},\qquad i = 1,2,\dots,n.
\end{gather*}
On the other hand, from \eqref{eq7.11} it follows that
\begin{gather*}
 \frac {\a_j(1-u_{ij})}{\a_j-u_{ij}}  = \frac {\a_k(1-u_{ik})}{\a_k-u_{ik}} \quad
\Rightarrow  \quad u_{ik}  = \frac {\a_k(1-\a_j)u_{ij}}{\a_j(1-\a_k)-(\a_j-\a_k)u_{ij}} .
\end{gather*}
So the equation to solve for $u_{ij}$, $j = 1,2,\dots,n$, for each $i$, is
\begin{gather*}
\frac {(1-\a_j)}{u_{ij}-\a_j} = \b_j + (1-\a_j) \sum_{k \ne j}^n \frac {\a_k\b_k}{\a_j(1-\a_k) - (\a_j-\a_k)u_{ij}} ,
\end{gather*}
since we must have $u_{ij} \ne 0$, which has $n$ roots.  Detailed analysis of these roots is not of immediate interest to us.

For $n = 2$,
\begin{gather*}
\frac {1-\a_1}{u_{11}\!-\a_1} = \b_1 + \frac {\a_2\b_2(1-\a_1)}{\a_1(1-\a_2)-(\a_1\!-\a_2)u_{11}}
= \frac {\a_1\b_1(1-\a_2)\! +\a_2\b_2(1-\a_2)\b_1(\a_1\!-\a_2)u_{11}}{\a_1(1-\a_2)-(\a_1\!-\a_2)u_{11}} ,\!\!
\end{gather*}
which leads to
\begin{gather*}
\b_1(\a_1-\a_2)(u_{11}-\a_1)^2 - (1-\a_1)(\a_1-\a_2+\a_1\b_1 + \a_2\b_2)(u_{11}-\a_1)
   + \a_1(1-\a_1)^2 = 0,
\end{gather*}
which is exactly the same as formula~(5.7) in \cite{HR2}.

Finally, it is clear from \eqref{eq7.9} that the eigenvalues are
\[
\l_{\bm} = \prod_{i=1}^n (1-\omega_i)^{m_i}.
\]

\subsection*{Acknowledgements}
The research of the f\/irst author was supported in part by the Applied Math. Sciences subprogram of the Of\/f\/ice of Energy Research, USDOE, under Contract DE-AC03-76SF00098.

\pdfbookmark[1]{References}{ref}
\LastPageEnding


\begin{thebibliography}{99}

\footnotesize\itemsep=0pt



\bibitem{AAR}
Andrews G., Askey R., Roy R.,
Special functions,
{\it Encyclopedia of Mathematics and its Applications}, Vol.~71, Cambridge University Press, Cambridge, 1999.


\bibitem{AK}
Aomoto K., Kita M.,
Hypergeometric funtions, Springer, Berlin, 1994 (in Japanese).

\bibitem{GI}
Geronimo J.S., Iliev P.,
Bispectrality of multivariable Racah--Wilson poynomials,
\href{http://dx.doi.org/10.1007/s00365-009-9045-3}{{\it Constr. Approx.}} {\bf 15} (2010), 417--457,
\href{http://arxiv.org/abs/0705.1469}{arXiv:0705.1469}.



\bibitem{G}
Gelfand I.M.,
General theory of hypergeometric functions,
{\it Soviet Math. Dokl.} {\bf 33} (1986), 573--577.

\bibitem{Gri}
Grif\/f\/iths R.C.,
Orthogonal polynomials on the multinomial distribution,
\href{http://dx.doi.org/10.1111/j.1467-842X.1971.tb01239.x}{{\it Austral.~J. Statist.}} {\bf 13} (1971),   27--35.

\bibitem{Gr}
Gr\"unbaum F.A.,
Block tridiagonal matrices and a beefed-up version of the Ehrenfest urn model,
in Modern Analysis and Applications, The Mark Krein Centenary Conference, Vol.~1,
Operator Theory and Related Topics,
\href{http://dx.doi.org/10.1007/978-3-7643-9919-1_15}{{\em Oper. Theory Adv. Appl.}}, Vol.~190, Birkh\"auser Verlag, Basel, 266--277.

\bibitem{GPT}
Gr\"unbaum F.A., Pacharoni I., Tirao J.A.,
Two stochastic models of a random walk in the $U(n)$-spherical duals of $U(n+1)$,
\href{http://arxiv.org/abs/1010.0720}{arXiv:1010.0720}.

\bibitem{HR1}
Hoare M.R., Rahman M.,
Cumultive Bernoulli trials and Krawtchouk processes,
\href{http://dx.doi.org/10.1016/0304-4149(84)90014-0}{{\it Stochastic Process. Appl.}} {\bf 16} (1983), 113--139.

\bibitem{HR2}
Hoare M.R., Rahman M.,
A probabilistic origin for a new class of bivariate polynomials,
\href{http://dx.doi.org/10.3842/SIGMA.2008.089}{{\it SIGMA}} {\bf 4} (2008), 089, 18~pages,
\href{http://arxiv.org/abs/0812.3879}{arXiv:0812.3879}.

\bibitem{HTF}
Erd\'elyi A., Magnus W., Oberhettinger F., Tricomi  F.G.,
Higher transcendental functions, Vols.~I,~II, III, McGraw-Hill Book Company, Inc., New York~-- Toronto~-- London, 1953, 1955.

\bibitem{IT}
Iliev P., Terwilliger P.,
The Rahman polynomials and the Lie algebra $sl_3(C)$,
\href{http://arxiv.org/abs/1006.5062}{arXiv:1006.5062}.

\bibitem{I}
Iliev P.,
A Lie theoretic interpretation of multivariate hypergeometric polynomials,
\href{http://arxiv.org/abs/1101.1683}{arXiv:1101.1683}.


\bibitem{M}
Mizukawa H.,
Zonal spherical functions on the complex ref\/lection groups and $(n+1,m+1)$-hypergeometric functions,
\href{http://dx.doi.org/10.1016/S0001-8708(03)00092-6}{{\it Adv. Math.}} {\bf 184} (2004),  1--17.

\bibitem{M1}
Mizukawa H.,
Orthogonality relations for multivariate Krawtchouck polynomials,
\href{http://dx.doi.org/10.3842/SIGMA.2011.017}{{\it SIGMA}} {\bf 7} (2011), 017, 5~pages,
\href{http://arxiv.org/abs/1009.1203}{arXiv:1009.1203}.

\bibitem{MT}
Mizukawa H., Tanaka H.,
$(n+1,m+1)$-hypergeometric functions associated to character algebras,
\href{http://dx.doi.org/10.1090/S0002-9939-04-07399-X}{{\it Proc. Amer. Math. Soc.}} {\bf 132} (2004), 2613--2618.


\bibitem{T}
Tirao J.A.,
The matrix-valued hypergeometric equation,
\href{http://dx.doi.org/10.1073/pnas.1337650100}{{\it Proc. Natl. Acad. Sci. USA}} {\bf 100} (2003), no.~14, 8138--8141.

\end{thebibliography}
\end{document}